\documentclass[11pt]{article}
\usepackage[T1]{fontenc}
\usepackage{latexsym,amssymb,amsmath,amsfonts,amsthm}
\makeatletter
\newcommand{\md}{\mbox{d}}
\newcommand{\be}{\begin{eqnarray}}
\newcommand{\ee}{\end{eqnarray}}
\newcommand{\mr}{\mathbb{R}}
\newcommand{\mc}{\mathcal}
\newcommand{\mb}{\mathbb}
\newcommand{\mx}{\mbox}

\newcommand{\tr}{\triangle}

\newcommand{\ds}{\displaystyle}

\newcommand{\ben}{\begin{eqnarray*}}
\newcommand{\enn}{\end{eqnarray*}}

\newcommand{\om}{\omega}

\newcommand{\al}{\alpha}
\newcommand{\la}{\lambda}
\newcommand{\La}{\Lambda}

\newtheorem{thm}{\textbf Theorem}[section]
\newtheorem{lem}{\textbf Lemma}[section]

\newtheorem{cor}{\textbf Corollary}[section]

\newtheorem{defin}{\textbf Definition}[section]

\begin{document}
\begin{titlepage}
\title{\bf Blow-up criterion of strong solutions to the Navier-Stokes equations in
Besov spaces with negative indices}
\author{Baoquan Yuan and Bo Zhang\\
         Institute of Applied Mathematics,\\
         Academy of Mathematics \& Systems Science\\
         Chinese Academy of Sciences\\
         Beijing 100080, P.R. China\\
          ({\sf bqyuan@hpu.edu.cn} (BY) and {\sf b.zhang@amt.ac.cn} (BZ))}
\date{}
\end{titlepage}
\maketitle

\begin{abstract}
A H\"older type inequality in Besov spaces is established and
applied to show that every strong solution $u(t,x)$ on (0,T) of the
Navier-Stokes equations can be continued beyond $t>T$ provided that
the vorticity $\omega(t,x)\in L^{\frac
2{2-\alpha}}(0,T;\dot{B}^{-\alpha}_{\infty,\infty}(\mr^3))\cap
L^{\frac2{1-\alpha}}(0,T;\dot{B}^{-1-\alpha}_{\infty,\infty}(\mr^3))$
for $0<\alpha<1$. \vskip0.1in

\noindent{\bf AMS Subject Classification 2000:} 35Q30.

\end{abstract}

\vspace{.2in} {\bf Key words:}\quad Navier-Stokes equations, Besov
spaces, blow-up criterion, regularity of weak solutions.



\section{Introduction}
\setcounter{equation}{0}

In this paper we consider the blow-up criterion of strong solutions to the Navier-Stokes equations
in $\mr^n$:
\begin{eqnarray}\label{1.1}
  \begin{cases}
   \partial_t u-\triangle u+(u\cdot\nabla)u+\nabla p=0, & x\in \mathbb{R}^n,\ t\in (0,T),\\
   \nabla\cdot u=0, & x\in \mathbb{R}^n ,\ t\in (0,T),\\
   u(0,x)=u_0(x), & x\in \mathbb{R}^n,
  \end{cases}
\end{eqnarray}
where $u=u(t,x)=(u^1(t,x),u^2(t,x),\cdots,u^n(t,x))$ and
$p=p(t,x)$ denote the unknown velocity vector and the unknown
pressure of the fluid at the point $(t,x)\in (0,T)\times\mr^n$,
respectively, and $u_0(x)=(u^1_0(x),u^2_0(x),\cdots,u^n_0(x))$ is
a given initial datum satisfying the divergence free constraint that $\nabla\cdot u_0(x)=0$,
and $n\ge 2$ is the space dimension.

There are many classical results on the local existence of smooth
solutions to the Navier-Stokes equations (\ref{1.1}).
For example, Fujita and Kato \cite{F-K} proved that for $u_0(x)\in H^s$ with
$s>\frac n2-1$ there exists a $T=T(\|u_0\|_{H^s})>0$ such that there exists
a unique solution $u(t,x)$ to the problem (\ref{1.1}) in the class
\begin{eqnarray}\label{cls}
u(t,x)\in C([0,T;H^s(\mr^n)))\cap C^1((0,T;H^s(\mr^n)))\cap
C((0,T;H^{s+2}(\mr^n))).
\end{eqnarray}
Later on, many authors established the well-posedness of the Navier-Stokes
equations (\ref{1.1}) in more general spaces \cite{B,K1,K2,KY,M-Y,T}.
It is interesting to ask whether or not the solution $u(t,x)$ blows up at time $t=T$ in some space-time
spaces. This question has been studied by many authors. For example,
Giga \cite{Y.G.} proved that, if $u(t,x)$ is a solution of the Navier-Stokes equations
(\ref{1.1}) in the class (\ref{cls}) and satisfies that
\begin{eqnarray}\label{Lpq}
u(t,x)\in L^q(0,T;L^p(\mr^n)),\qquad \frac2q+\frac np=1,\quad n<p\le\infty,\ 2\le q<\infty,
\end{eqnarray}
then $u$ can be extended, as a solution in the class (\ref{cls}),
onto $(0,T')$ for some $T'>T$. On the other hand, the space-time
Lebesgue space $L^q(0,T;L^p(\mathbb{R}^n))$ plays a significant
role in the study of regularity of weak solutions. For instance,
Leray-Hopf weak solution $u(t,x)$ becomes a unique strong solution
of the Navier-Stokes equations (\ref{1.1}) on $(0,T]$ if it satisfies the condition
(\ref{Lpq}). For regularity of weak solutions, see \cite{SR,F-J-R,K-O-T}.

From the hydrodynamical point of view, we want to study the regularity of weak solutions
by means of the vorticity $\omega(t,x)=\mbox{curl} u(t,x)=\nabla\times u(t,x)$.
Beir$\tilde{a}$o da Veiga proved that, if the vorticity $\omega(t,x)$ of the Leray-Hopf weak solution
$u(t,x)$ belongs to $L^q(0,T;L^p(\mathbb{R}^n))$ with $\ds\frac2q+\frac np=2$
for $1<q<\infty$ and $n/2<p<\infty$, then $u(t,x)$ becomes a strong solution on $(0,T)$ (see \cite{B-V}).
In the end-point case $p=\infty$, Kozono and Taniuchi \cite{K-T} established the regularity
of weak solutions under the condition
\ben
\omega(t,x)\in L^1(0,T;BMO),
\enn
where BMO stands for the space of bounded mean oscillations.
Based on a logarithmic Sobolev inequality in Besov spaces, Kozono, Ogawa and Taniuchi \cite{K-O-T}
extended the above condition to the condition
\ben
\omega(t,x)\in L^q(0,T;\dot{B}^0_{p,\infty}(\mr^n))\ \mx{with }
\frac 2q+\frac np=2,\ \frac n2<p\le \infty,
\enn
where $\dot{B}^s_{p,q}(\mr^n)$ stands for the homogeneous Besov space.
Kozono and Shimada \cite{K-S} established the continuation principle of the smooth solution
$u(t,x)$ on $(0,T)$ under the condition
\be\label{1.6}
u(t,x)\in L^{\frac2{1-\alpha}}(0,T;\dot{F}^{-\alpha}_{\infty,\infty}(\mr^n)),\qquad
0<\alpha<1,
\ee
where $\dot{F}^s_{p,q}(\mr^n)$ stands for the homogeneous Triebel-Lizorkin space.

From the scaling invariant point of view, it is important that the solution $u(t,x)$
in the space $L^q(0,T;L^p(\mr^n))$ is scaling invariant, that is,
$\|u_\lambda(t,x)\|_{L^q(0,T;L^p)}=\|u(t,x)\|_{L^q(0,T;L^p)}$ for
$2/q+n/p=1$, $2\le q<\infty$ and $n<p\le \infty,$ where $u_\la(t,x)=\la u(\la^2 t,\la x)$
with $\lambda>0.$
Similarly, the vorticity $\omega(t,x)$ of a solution is also scaling invariant, that is,
$\|\omega_\lambda(t,x)\|_{L^q(0,T;L^p)}=\|\omega(t,x)\|_{L^q(0,T;L^p)}$
for $2/q+n/p=2$, $1\le q<\infty$ and $n/2<p\le\infty$, where
$\omega_\lambda(t,x)=\lambda^2\omega(\lambda^2 t,\lambda x)$ with $\lambda>0$.
So the regularity criterion of weak solutions in terms of the solution $u(t,x)$ and the
vorticity $\omega(t,x)$ is very natural by the scaling invariance standard.

It is natural to ask if the continuation principle remains true for the smooth solution $u(t,x)$ on $(0,T)$
under a condition for the vorticity $\omega(t,x)$ similar to (\ref{1.6}).
In this paper we shall prove that such a continuation criterion can be derived
for the smooth solution $u(t,x)$ on $(0,T)$ under the assumption that
\be\label{1.7}
\om(t,x)\in L^{\frac2{2-\al}}(0,T;\dot{B}^{-\al}_{\infty,\infty}(\mr^3))
\cap L^{\frac2{1-\al}}(0,T;\dot{B}^{-1-\al}_{\infty,\infty}(\mr^3))
\ee
for $0<\al<1$.  Note that the space-time spaces
$L^{\frac2{1-\al}}(0,T;\dot{F}^{-\al}_{\infty,\infty}(\mr^n))$
and $L^{\frac2{2-\al}}(0,T;\dot{B}^{-\al}_{\infty,\infty}(\mr^n))
\cap L^{\frac2{1-\al}}(0,T;\dot{B}^{-1-\al}_{\infty,\infty}(\mr^n))$
are scaling invariant.

To establish the continuation principle, we need a bilinear estimate
of the H\"older type in the homogeneous Besov spaces
$\dot{B}^s_{p,q}(\mr^n)$ for $s>0$, $1\le p,\ q\le \infty$:
\be\label{1.8}
\|f_1\cdot f_2;\dot{B}^s_{p,q}\|\le C(\|f_1;\dot{B}^{s+\al}_{p_1,q_1}\|\|f_2;\dot{B}^{-\al}_{p_2,q_2}\|
+\|f_1;\dot{B}^{-\beta}_{p_3,q_3}\|\|f_2;\dot{B}^{s+\beta}_{p_4,q_4}\|)
\ee
for any $f_1\in\dot{B}^{s+\al}_{p_1,q_1}\cap\dot{B}^{-\beta}_{p_3,q_3}$,
$f_2\in \dot{B}^{s+\beta}_{p_4,q_4}\cap\dot{B}^{-\al}_{p_2,q_2}$,
where $\alpha,\ \beta>0$ and $s>0$, $1/p=1/p_1+1/p_2=1/p_3+1/p_4$
and $1/q=1/q_1+1/q_2=1/q_3+1/q_4$.

The paper is organized as follows. In section 2, we first recall the definition of
Besov spaces and then prove the H\"older type inequality (\ref{1.8}).
Section 3 is devoted to the proof of the main theorem on the continuation principle
under the condition (\ref{1.7}) for the strong solution $u(t,x)$.

%


\section{Preliminaries}
\setcounter{equation}{0}

We first introduce the Littlewood-Paley decomposition and the definition of Besov spaces.
For $f\in\mc{S}(\mr^n),$ the Schwartz class of rapidly decreasing functions,
define the Fourier transform
\ben
\hat{f}(\xi)=\mc{F}f(\xi)=(2\pi)^{-n/2}\int_{\mr^n}\mx{e}^{-ix\cdot\xi}f(x)\md x
\enn
and the inverse Fourier transform:
\ben
\check{f}(x)=\mc{F}^{-1}f(x)=(2\pi)^{-n/2}\int_{\mr^n}\mx{e}^{ix\cdot\xi}f(\xi)\md\xi.
\enn
Choose a non-negative radial function $\chi(\xi)\in C^\infty_0(\mr^n)$ such that
$0\le \chi(\xi)\le 1$ and
\ben
\chi(\xi)=\begin{cases}
1, &\mx{ for } |\xi|\le \frac34,\\
0, &\mx{ for } |\xi|>\frac 43,
\end{cases}
\enn
and let $\hat{\varphi}(\xi)=\chi(\xi/2)-\chi(\xi)$,
$\chi_j(\xi)=\chi(\frac{\xi}{2^j})$ and
$\hat{\varphi}_j(\xi)=\hat{\varphi}(\frac \xi{2^j})$ for $j\in\mb{Z}$.
Write
\ben
h(x)&=&\mc{F}^{-1}\chi(x),\\
h_j(x)&=&2^{nj}h(2^jx),\\
\varphi_j(x)&=&2^{nj}\varphi(2^jx).
\enn
Define the Littlewood-Paley projection operators $S_j$ and
$\tr_j$, respectively, as
\ben
S_ju(x)&=&h_j*u(x)\ \mx{ for } j\in \mb{Z},\\
\tr_ju(x)&=&\varphi_j*u(x)=S_{j+1}u(x)-S_ju(x) \ \mx{ for } j\in\mb{Z}.
\enn
Formally, $\tr_j$ is a frequency projection to the annulus $|\xi|\sim 2^j$,
whilst $S_j$ is a frequency projection to the ball $|\xi|\lesssim 2^j$ for $j\in \mb{Z}$.
For any $u(x)\in L^2(\mr^n)$ we have the Littlewood-Paley decomposition
\ben
u(x)&=&h*u(x)+\sum_{j\ge 0}\varphi_j*u(x),\\
u(x)&=&\sum^\infty_{j=-\infty}\varphi*u(x).
\enn
Clearly,
\ben
\mx{supp}\,\chi(\xi)\cap\mx{supp}\,\hat{\varphi}_j(\xi)&=&\emptyset \mx{ for } j\ge 1,\\
\mx{supp}\,\hat{\varphi}_j(\xi)\cap\mx{supp}\,\hat{\varphi}_{j'}(\xi)&=&\emptyset,
\mx{ for } |j-j'|\ge 2.
\enn

We now recall the definition of Besov spaces.
Let $s\in \mr$ and $1\le p,\ q\le \infty$, the Besov space $B^s_{p,q}(\mr^n)$ (denote by
$B^s_{p,q}$) is defined by
\ben
B^s_{p,q}=\{f\in\mc{S}(\mr^n):\,\|f\|_{B^s_{p,q}}<\infty \}
\enn
with the norm
\ben
\|f\|_{B^s_{p,q}}=(\|h*f\|^q_p+\sum_{j\ge 0}2^{jsq}\|\varphi_j*f\|^q_p)^{1/q}.
\enn
The homogeneous Besov space $\dot{B}^s_{p,q}$ is defined by the dyadic decomposition as
\ben
\dot{B}^s_{p,q}=\{f\in \mathcal{Z}'(\mathbb{R}^n):\,\|f\|_{\dot{B}^s_{p,q}}<\infty\}
\enn
with the norm
\ben
\|f\|_{\dot{B}^s_{p,q}}=(\sum^\infty_{j=-\infty}2^{jsq}\|\varphi_j*f\|^q_p)^{1/q},
\enn
where $\mathcal{Z}'(\mathbb{R}^n)$ denotes the dual space of
$$\mathcal{Z}(\mathbb{R}^n)=\{f\in\mathcal{S}(\mathbb{R}^n):\,D^\al\hat{f}(0)=0,\,
\mbox{for any multi-index}\,\al\in \mathbb{N}^n\}$$
and can be identified by the quotient space $\mathcal{S}'/\mathcal{P}$ with
the set $\mathcal{P}$ of polynomial functions.
For details see \cite{B-L}, \cite{M1} or \cite{Tr}.

For convenience we recall the definition of Bony's para-product formula which gives
the decomposition of the product $f\cdot g$ of two functions $f(x)$ and $g(x)$.
%

\begin{defin}
The para-product of two functions $f$ and $g$ is defined by
\ben
T_gf=\sum_{i\le j-2}\triangle_ig\triangle_jf=\sum_{j\in
\mathbb{Z}}S_{j-1}g\triangle_jf.
\enn
The remainder of the para-product is defined by
\ben
R(f,g)=\sum_{|i-j|\le1}\triangle_ig\triangle_jf.
\enn
Then Bony's para-product formula reads
\begin{eqnarray}\label{2.18}
f\cdot g=T_gf+T_fg+R(f,g).
\end{eqnarray}
\end{defin}

We now have the bilinear estimate of the H\"older type
in the Besov space $\dot{B}^s_{p,q}$.

\begin{lem}\label{lem2.1}
Let $1\le p,\ q\le\infty$, $s>0$, $\al>0$ and $\beta>0$.
Choose $1\le p_i,\ q_i\le\infty$ ($i=1,\ 2,\ 3,\ 4$) so that
\ben
\frac1p=\frac1{p_1}+\frac1{p_2}=\frac1{p_3}+\frac1{p_4},\\
\frac1q=\frac1{q_1}+\frac1{q_2}=\frac1q_3+\frac1{q_4}.
\enn
Then there exists a constant $C$ such that $f_1\cdot f_2\in\dot{B}^s_{p,q}(\mr^n)$ and
\ben
\|f_1\cdot f_2;\dot{B}^s_{p,q}\|\le C(\|f_1;\dot{B}^{s+\al}_{p_1,q_1}\|\|f_2;\dot{B}^{-\al}_{p_2,q_2}\|
+\|f_1;\dot{B}^{-\beta}_{p_3,q_3}\|\|f_2;\dot{B}^{s+\beta}_{p_4,q_4}\|)
\enn
for any $f_1\in\dot{B}^{s+\al}_{p_1,q_1}\cap\dot{B}^{-\beta}_{p_3,q_3}$,
$f_2\in\dot{B}^{s+\beta}_{p_4,q_4}\cap\dot{B}^{-\alpha}_{p_2,q_2}$.
\end{lem}

\begin{proof}
By Bony's para-product decomposition (\ref{2.18}) one has
\be\nonumber
f_1\cdot f_2&=&\sum_{k=-\infty}^{\infty}\tr_kf_1S_{k-1}f_2+\sum_{k=-\infty}^{\infty}
\tr_kf_2S_{k-1}f_1+\sum_{|i-j|\le1}\tr_if_1\tr_jf_2\\ \label{2.22}
&\triangleq& I_1+I_2+I_3.
\ee
We first compute the Besov norm $\|I_1\|_{\dot{B}^s_{p,q}}$ of $I_1$.
Since supp$[\mc{F}(S_{k-1}f_2)]\subseteqq\{|\xi|\le\frac23 2^k\}$ and
supp$[\mc{F}(\tr_kf_1)]\subseteqq \{\frac342^k\le|\xi|\le\frac23 2^{k+2}\}$,
then we have
\ben
\mx{supp}[\mc{F}(\tr_kf_1S_{k-1}f_2)]\subseteqq\Big\{\frac13 2^{k-2}\le|\xi|\le\frac532^{k+1}\Big\}.
\enn
Thus by the definition of the Besov norm we can derive that
\be\nonumber
\|I_1\|_{\dot{B}^s_{p,q}}&=&(\sum_{j\in\mb{Z}}2^{sjq}\|\tr_j\sum_{k\in\mb{Z}}
\tr_kf_1S_{k-1}g\|^q_p)^{1/q}\\ \nonumber
&\le& C(\sum_{j\in\mb{Z}}2^{sjq}\|\tr_j\sum_{k=j-3}^{j+5}\tr_kf_1S_{k-1}f_2\|^q_p)^{1/q}\\ \label{2.24}
&\le& C\bigg(\sum_{j\in\mb{Z}}2^{sjq}\Big(\sum_{l=-3}^{5}\|\tr_{l+j}f_1\|_{p_1}
\|S_{j+l-1}f_2\|_{p_2}\Big)^q\bigg)^{1/q},
\ee
where use has been made of the H\"older inequality and the fact that
$\|\tr_jf\|_p\le C\|f\|_p$ for any $j\in \mb{Z}$.
By the Minkowski and H\"older inequalities the last term on the right-hand side of (\ref{2.24})
can be bounded above by
\ben
&& C\sum_{l=-3}^5\Big(\sum_{j\in\mb{Z}}2^{sjq}\|\tr_{l+j}f_1\|^q_{p_1}\|S_{l+j-1}f_2\|^q_{p_2}\Big)^{1/q}\\
&\le& C\sum_{l=-3}^5 2^{-\al-sl}\Big(\sum_{j\in\mb{Z}}2^{(s+\al)(j+l)}
\|\tr_{l+j}f_1\|^q_{p_1}2^{-\al(l+j-1)q}\|S_{l+j-1}f_2\|^q_{p_2}\Big)^{1/q}\\
&\le& C\|f_1\|_{\dot{B}^{s+\al}_{p_1,q_1}}\Big(\sum_{j\in\mb{Z}}2^{-\al(l+j-1)q_2}
\|S_{l+j-1}f_2\|^{q_2}_{p_2}\Big)^{1/q_2}.
\enn
Noting the equivalence of the Besov norms
$$\Big(\sum_{j\in\mb{Z}}(2^{sj}\|\tr_jf\|_p)^q\Big)^{1/q}
\thicksim\Big(\sum_{j\in\mb{Z}}(2^{sj}\|S_jf\|_p)^q\Big)^{1/q}$$
for $s<0,$ we arrive at
\be\label{2.26}
\|I_1\|_{\dot{B}^s_{p,q}}\le C\|f_1\|_{\dot{B}^{s+\al}_{p_1,q_1}}\|f_2\|_{\dot{B}^{-\al}_{p_2,q_2}}.
\ee
If either $p$ or $q$ equals to infinite, the above argument can be modified accordingly to get
the estimate (\ref{2.26}). For $I_2$ arguing similarly as above gives
\be\label{2.27}
\|I_2\|_{\dot{B}^s_{p,q}}\le C\|f_2\|_{\dot{B}^{s+\beta}_{p_3,q_3}}
\|f_1\|_{\dot{B}^{-\beta}_{p_4,q_4}}.
\ee
We now consider $I_3$.
Since supp$[\mc{F}(\tr_kf_1\tr_{k+l}f_2)]\subseteqq \{|\xi|\le\frac83 2^k(1+2^l)\}$ and
supp$[\mc{F}(\tr_jf)]\subseteqq\{\frac34 2^j\le |\xi|\le \frac83 2^j\}$, it follows by the definition
of Besov spaces that
\ben\nonumber
\|I_3\|_{\dot{B}^s_{p,q}}&=& \bigg(\sum_{j\in\mb{Z}}2^{sjq}\Big\|\tr_j
\sum_{|k-l|\le1}\tr_kf_1\tr_lf_2\Big\|^q_p\bigg)^{1/q}\\ \nonumber
&\le& C\bigg(\sum_{j\in\mb{Z}}2^{sjq}\Big\|\tr_j\sum_{k\ge j-4}
\sum_{l=-1}^1\tr_kf_1\tr_{k+l}f_2\Big\|^q_p\bigg)^{1/q}\\ \label{2.28}
&\le& C\bigg(\sum_{j\in \mb{Z}}2^{sjq}\Big(\sum_{r\ge-4}
\sum_{l=-1}^1\|\tr_{j+r}f_1\|_{p_1}\|\tr_{r+j+l}f_2\|_{p_2}\Big)^q\bigg)^{1/q},
\enn
where use has been made of the Minkowski and H\"older inequalities.
Using the Minkowski and H\"older inequalities again the above inequality
can be evaluated continually as
\be\nonumber
&&\|I_3\|_{\dot{B}^s_{p,q}}\le C\sum_{r\ge -4}\sum_{l=-1}^1\bigg(\sum_{j\in\mb{Z}}2^{sjq}
\|\tr_{j+r}f_1\|^q_{p_1}\|\tr_{r+j+l}f_2\|^q_{p_2}\bigg)^{1/q}\\ \nonumber
&&\le C\sum_{r\ge-4}\sum_{l=-1}^1\bigg(2^{(-sr+\alpha l)q}
\sum_{j\in\mb{Z}}2^{(s+\alpha)(j+r)q}\|\tr_{j+r}f_1\|^q_{p_1}2^{-\alpha(j+r+l)q}
\|\tr_{j+r+l}f_2\|^q_{p_2}\bigg)^{1/q}\\ \nonumber
&&\le C\sum_{r\ge -4}\sum_{l=-1}^1 2^{-sr+\al l}\bigg(\sum_{j\in\mb{Z}}2^{(s+\alpha)jq_1}
\|\tr_jf_1\|^{q_1}_{p_1}\bigg)^{1/{q_1}}\bigg(\sum_{j\in\mb{Z}}2^{-\al jq_2}
\|\tr_jf_2\|^{q_2}_{p_2}\bigg)^{1/{q_2}}\\ \label{2.29}
&&\le C\|f_1\|_{\dot{B}^{s+\alpha}_{p_1,q_1}}\|f_2\|_{\dot{B}^{-\alpha}_{p_2,q_2}}.
\ee
Collecting the estimates (\ref{2.22}), (\ref{2.26}), (\ref{2.27}) and
(\ref{2.29}) completes the proof of Lemma \ref{lem2.1}.
\end{proof}


\section{The blow-up criterion}
\setcounter{equation}{0}

In this section we take $n=3$ as an example to present the blow-up criterion.
We first introduce some notations and function spaces.
Denote by $C^\infty_{0,\sigma}(\mr^3)$ the set of all $C^\infty$ vector functions
$f(x)=(f_1(x),f_2(x),f_3(x))$ with compact support satisfying that $\nabla\cdot f(x)=0$.
$L^r_\sigma(\mr^3)$ is the closure of $C^\infty_{0,\sigma}(\mr^3)$-functions with respect to
the $L^r$-norm $\|\cdot\|_r$ for $1\le r\le \infty$.
$H^s_\sigma(\mr^3)$ denotes the closure of $C^\infty_{0,\sigma}(\mr^3)$-functions with respect to the
$H^s$-norm $\|f\|_{H^s}=\|(1-\tr)^{s/2}f\|_2$ for $s\ge 0$.
We are now ready to state and prove the continuation principle of strong solutions to
the Navier-Stokes equations.

\begin{thm}\label{thm3.1}
Let $0<\alpha<1$, $T>0$ and $u_0\in H^1_\sigma(\mr^3)$.
Assume that
\ben
u(t,x)\in C([0,T);H^1_\sigma(\mr^3))\cap C^1((0,T);H^1_\sigma(\mr^3))
\cap C((0,T);H^3_\sigma(\mr^3))
\enn
is a strong solution to the Navier-Stokes equations (\ref{1.1}).
If $u(t,x)$ satisfies
\be\label{3.2}
\int^T_0\bigg(\|\omega(\tau,\cdot)\|^{\frac2{2-\al}}_{\dot{B}^{-\al}_{\infty,\infty}}
+\|\omega(\tau,\cdot)\|^{\frac2{1-\al}}_{\dot{B}^{-1-\al}_{\infty,\infty}}\bigg)\md\tau<\infty,
\ee
then $u(t,x)$ can be continually extended to the interval
$(0,T')$ for some $T'>T$, where $\omega(t,x)=\nabla\times u(t,x)$ is the vorticity.
\end{thm}

An immediate corollary is as follows.

\begin{cor}\label{cor3.1}
Let $0<\alpha<1$ and let $u(t,x)$ be a strong solution to the
Navier-Stokes equations (\ref{1.1}) satisfying that
\ben
u(t,x)\in C([0,T);H^1_\sigma(\mr^3))\cap C^1((0,T);H^1_\sigma(\mr^3))\cap C((0,T);H^3_\sigma(\mr^3)).
\enn
If $T$ is the maximal existence time, then
\be\label{3.4}
\int^T_0\bigg(\|\omega(\tau,\cdot)\|^{\frac2{2-\al}}_{\dot{B}^{-\al}_{\infty,\infty}}
+\|\omega(\tau,\cdot)\|^{\frac2{1-\al}}_{\dot{B}^{-1-\al}_{\infty,\infty}}\bigg)\md\tau=\infty.
\ee
\end{cor}

{\it Proof of Theorem \ref{thm3.1}.}
We begin with the Navier-Stokes equations in the vorticity form
\begin{eqnarray}\label{NSV}
 \left\{
  \begin{array}{ll}
   \partial_t \omega-\triangle \omega+(u\cdot\nabla)\omega-(\omega\cdot\nabla)u=0,
   \quad x\in \mathbb{R}^3,\ t\in (0,T),\\[0.2cm]
    \omega(x,0)=\omega_0(x),\ x\in \mathbb{R}^3, \\[0.2cm]
     \mathrm {div} \omega_0(x)=0,\ x\in \mathbb{R}^3.
  \end{array}
 \right.
\end{eqnarray}
Multiplying the first equation of (\ref{NSV}) by $\omega(t,x)$ and
integrating by parts, it follows that \be\label{3.6}
\frac12\frac{\md}{\md
t}\|\omega(t,\cdot)\|_2^2+\|\nabla\omega(t,\cdot)\|^2_2
=(\omega\cdot\nabla u,\omega), \ee where use has been made of the
fact that $(u\cdot\nabla\omega,\omega)=0$. Since $(\omega\cdot\nabla
u,\omega)=(\nabla(\omega\otimes u),\omega)$, we have by the
Cauchy-Schwarz inequality that \be\nonumber |(\omega\cdot\nabla
u,\omega)|=|(\nabla(\om\otimes u),\om)|
&=&|(\La^{-\al}\nabla(\om\otimes u),\La^{\al}\om)|\\ \label{3.7}
&\le& \|\omega\otimes
u\|_{\dot{B}^{1-\alpha}_{2,2}}\|\|\Lambda^{\alpha}\omega\|_2, \ee
where $\Lambda=(-\tr)^{1/2}$.

We now estimate $\|\omega\otimes u\|_{\dot{B}^{1-\alpha}_{2,2}}$. By
Lemma \ref{lem2.1} with $\alpha=\beta$, $p_1=q_1=2$,
$p_2=q_2=\infty$, $p_3=q_3=\infty$ and $p_4=q_4=2$ it follows that
\be\nonumber \|\om\otimes u\|_{\dot{B}^{1-\al}_{2,2}}&\le&
C(\|\om\|_{\dot{B}^1_{2,2}} \|u\|_{\dot{B}^{-\al}_{\infty,\infty}}
+\|\omega\|_{\dot{B}^{-\alpha}_{\infty,\infty}}\|u\|_{\dot{B}^1_{2,2}})\\
\label{3.8} &\le& C(\|\nabla\omega\|_2\|\nabla
u\|_{\dot{B}^{-1-\alpha}_{\infty,\infty}}
+\|\omega\|_{\dot{B}^{-\alpha}_{\infty,\infty}}\|\nabla
u\|_{\dot{B}^0_{2,2}}). \ee Noting the $L^p$ boundedness of singular
integral operators for $1<p<\infty$: \be\label{3.9} \|\nabla
u(t,\cdot)\|_p\le C\|\omega(t,\cdot)\|_p, \ee and by (\ref{3.7}),
(\ref{3.8}) and the Gagliardo-Nirenberg inequality we obtain that
\ben\nonumber |(\omega\cdot\nabla u,\omega)|&\le&
C(\|\nabla\omega\|_2\|
\omega\|_{\dot{B}^{-1-\alpha}_{\infty,\infty}}
+\|\om\|_{\dot{B}^{-\al}_{\infty,\infty}}\|\om\|_2)\|\La^{\al}\om\|_2\\
\nonumber &\le&
C(\|\nabla\omega\|_2\|\omega\|_{\dot{B}^{-1-\alpha}_{\infty,\infty}}
+\|\omega\|_{\dot{B}^{-\alpha}_{\infty,\infty}}\|\omega\|_2)
\|\omega\|_2^{1-\alpha}\|\nabla\omega\|_2^\alpha\\ \label{3.10}
&\le& C(\|\nabla\omega\|^{1+\alpha}_2\|\omega\|^{1-\alpha}_2
\|\omega\|_{\dot{B}^{-1-\alpha}_{\infty,\infty}}+\|\nabla\omega\|^\alpha_2
\|\omega\|^{2-\alpha}_2\|\omega\|_{\dot{B}^{-\alpha}_{\infty,\infty}}).
\enn By Young's inequality the above inequality can be estimated
continually as \be\nonumber |(\omega\cdot\nabla
u,\omega)|&\le&\frac14\|\nabla\omega\|^2_2
+C\|\omega\|^{\frac2{1-\alpha}}_{\dot{B}^{-1-\alpha}_{\infty,\infty}}\|\omega\|^2_2
+\frac14\|\nabla\om\|^2_2+\|\om\|^2_2\|\om\|^{\frac2{2-\al}}_{\dot{B}^{-\al}_{\infty,\infty}}\\
\label{3.11}
&\le&\frac12\|\nabla\om\|^2_2+C\|\om\|^2_2\Big(\|\om\|^{\frac2{1-\al}}_{\dot{B}^{-1-\al}_{\infty,\infty}}
+\|\omega\|^{\frac2{2-\alpha}}_{\dot{B}^{-\alpha}_{\infty,\infty}}\Big).
\ee Combining (\ref{3.11}) with (\ref{3.6}) we arrive at
\be\label{3.12} \frac12\frac{\md}{\md
t}\|\omega(t,\cdot)\|_2^2+\|\nabla\omega(t,\cdot)\|^2_2 \le
C\|\omega\|^2_2\Big(\|\omega\|^{\frac2{1-\alpha}}_{\dot{B}^{-1-\alpha}_{\infty,\infty}}
+\|\omega\|^{\frac2{2-\alpha}}_{\dot{B}^{-\alpha}_{\infty,\infty}}\Big).
\ee Applying the Gronwnall inequality to (\ref{3.12}) yields \ben
\|\omega(t\|_2\le\|\omega_0\|_2\exp\bigg\{-C\int^T_0
\Big(\|\omega(\tau)\|^{\frac2{1-\alpha}}_{\dot{B}^{-1-\alpha}_{\infty,\infty}}
+\|\omega(\tau)\|^{\frac2{2-\alpha}}_{\dot{B}^{-\alpha}_{\infty,\infty}}\Big)\md\tau\bigg\}.
\enn This, together with (\ref{3.9}) with $p=2$, implies that
\be\label{3.14} \|\nabla
u(t)\|_2\le\|\omega_0\|_2\exp\bigg\{-C\int^T_0
\Big(\|\omega(\tau)\|^{\frac2{1-\alpha}}_{\dot{B}^{-1-\alpha}_{\infty,\infty}}
+\|\omega(\tau)\|^{\frac2{2-\alpha}}_{\dot{B}^{-\alpha}_{\infty,\infty}}\Big)\md\tau\bigg\}.
\ee On the other hand, the strong solution $u(t,x)$ satisfies the
energy identity: \be\label{3.15} \|u(t)\|^2_2+2\int^t_0\|\nabla
u(\tau)\|^2_2\md\tau=\|u_0\|^2_2 \ee for $0\le t\le T$. Thus by
(\ref{3.14}) and (\ref{3.15}) we have \ben\nonumber
\|u(t)\|_{H^1}\le\|u_0\|_2+\|\omega_0\|_2\exp\bigg\{-C\int^T_0
\Big(\|\omega(\tau)\|^{\frac2{1-\alpha}}_{\dot{B}^{-1-\alpha}_{\infty,\infty}}
+\|\omega(\tau)\|^{\frac2{2-\alpha}}_{\dot{B}^{-\alpha}_{\infty,\infty}}\Big)\md\tau\bigg\}
<\infty, \enn provided that $u$ satisfies the condition
(\ref{3.2}). This completes the proof of Theorem \ref{thm3.1}.


\section*{Acknowledgements}

The research of B Yuan was partially supported by the National
Natural Science Foundation of China No. 10571016, and Natural
Science Foundation of Henan Province No. 0611055500. The research of
B Zhang was supported by the Chinese Academy of Sciences through the
Hundred Talents Program.


\end{document}